\documentclass[11pt]{article}

\usepackage{amssymb}
\usepackage{amscd}
\usepackage{amsmath}
\usepackage{amsfonts}
\usepackage{theorem}
\usepackage{mathrsfs}
\usepackage{tikz}
\usepackage{graphicx}
\usepackage{float}
\usepackage{color}
\definecolor{red}{rgb}{.7,0,0}
\definecolor{blue}{rgb}{0,0,1}

{\theorembodyfont{\slshape}
\newtheorem{theorem}{Theorem}
\newtheorem{proposition}{Proposition}
\newtheorem{lemma}{Lemma}

}
{\theorembodyfont{\rmfamily}

}
\def\proof{{\noindent\sc Proof. \quad}}
\newcommand{\proofof}[1]{{\noindent\sc Proof of #1. \quad}}
\def\eproof{{\mbox{}\hfill\qed}\medskip}
\newcommand\qed{{\unskip\nobreak\hfil\penalty50\hskip2em\vadjust{}
\nobreak\hfil$\Box$\parfillskip=0pt\finalhyphendemerits=0\par}}

\def\N{{\mathbb{N}}}
\def\C{{\mathbb{C}}}

\renewcommand{\P}{\mathbb{P}}
\def\E{\mathop{\mathbb{E}}}
\def\Avg{\mathop{\mathsf{Avg}}}
\def\std{\mathsf{std}}

\def\Cnn{{\C^{n\times n}}}

\def\cond{\mathsf{cond}}
\def\Cw{\mathsf{Cw}}

\def\condmin{\mathsf{cond}_{\min}}
\def\condmax{\mathsf{cond}_{\max}}

\def\Id{\mathsf{Id}}

\renewcommand{\tilde}{\widetilde}

\def\z{\zeta}

\def\Cnn{\C^{n\times n}}

\newcommand{\algoritmo}{\begin{minipage}{0.87\hsize}\linea}
\newcommand{\falgoritmo}{\linea\end{minipage}\bigskip}
\newcommand{\linea}{\vspace*{-5pt}\hrule\vspace*{5pt}}

\medskip

\begin{document}

\bibliographystyle{plain}

\makeatletter


\def\JACM{Journal of the ACM}
\def\CACM{Communications of the ACM}
\def\ICALP{International Colloquium on Automata, Languages
            and Programming}
\def\STOC{annual ACM Symp. on the Theory
          of Computing}
\def\FOCS{annual IEEE Symp. on Foundations of Computer Science}
\def\SIAM{SIAM Journal on Computing}
\def\SIOPT{SIAM Journal on Optimization}
\def\MOR{Math. Oper. Res.}
\def\BSMF{Bulletin de la Soci\'et\'e Ma\-th\'e\-ma\-tique de France}
\def\CRAS{C. R. Acad. Sci. Paris}
\def\IPL{Information Processing Letters}
\def\TCS{Theoretical Computer Science}
\def\BAMS{Bulletin of the Amer. Math. Soc.}
\def\TAMS{Transactions of the Amer. Math. Soc.}
\def\PAMS{Proceedings of the Amer. Math. Soc.}
\def\JAMS{Journal of the Amer. Math. Soc.}
\def\LNM{Lect. Notes in Math.}
\def\LNCS{Lect. Notes in Comp. Sci.}
\def\JSL{Journal for Symbolic Logic}
\def\JSC{Journal of Symbolic Computation}
\def\JCSS{J. Comput. System Sci.}
\def\JoC{J. of Complexity}
\def\MP{Math. Program.}
\sloppy

\begin{title}
{{\bf On the condition of characteristic polynomials}}
\end{title}

\author{Peter B\"urgisser\thanks{
Institute of Mathematics, Technische Universit\"at Berlin,
pbuerg@math.tu-berlin.de.
Partially supported by DFG grant BU 1371/2-2.}
\and
Felipe Cucker\thanks{Department of Mathematics,
City University of Hong Kong, macucker@cityu.edu.hk.
Partially supported by a GRF grant from the Research
Grants Council of Hong Kong (project number CityU~100813).}
\and
Elisa Rocha Cardozo\thanks{Centro de Matem\'atica,
Universidad de la Rep\'ublica, Uruguay, elisa@cmat.edu.uy. Partially supported by PEDECIBA.}
}
\date{\today}

\date{}

\makeatletter
\maketitle
\makeatother

\begin{quote}{\small
{\bf Abstract.}\quad
We prove that the expectation of the logarithm of the condition
number of each of the zeros of the characteristic polynomial
of a complex standard Gaussian matrix is $\Omega(n)$.
This may provide an explanation for the common wisdom in
numerical linear algebra that advises against computing
eigenvalues via root-finding for characteristic polynomials.
}
\end{quote}

\smallskip

\noindent{\bf AMS subject classifications:} 65F35, 65Y20, 68F15

\smallskip

\noindent{\bf Key words:} eigenvalues, characteristic polynomial,
condition, random matrices

\thispagestyle{empty}

\section{Introduction}\label{sec:intro}

Common wisdom in numerical linear algebra advices against
computing eigenvalues via root-finding for characteristic
polynomials. Thus, the chapter on eigenvalues
in Datta's textbook~\cite[p.~372]{Datta:95}
begins with the following lines:

\begin{quote}{\small
Because the eigenvalues of a matrix $A$ are the zeros of the
characteristic polynomial $\det(A-\lambda \Id)$, one would naively
think of computing the eigenvalues of~$A$ by finding its
characteristic polynomial and then computing its zeros by a
standard root-finding method. Unfortunately, this is not a
practical approach. }
\end{quote}
Unfortunately as well, there is no explanation on why this is not a
practical approach. Complexity not being an issue (there is a vast
literature on efficient algorithms for univariate polynomials
root finding~\cite{Pan97,Schoen82}),
it appears that numerical stability is.

According to Wilkinson~\cite[p.~13]{Wilkinson:84},
``almost all of the algorithms developed before the 1950's for
dealing with the unsymmetric eigenvalue problem [\dots] were
based on some device for computing the explicit polynomial
equation.'' But in the early fifties this practice was put into
question as some examples (notably the ``perfidious polynomial'')
showed an unexpectedly high sensitivity to root finding (i.e.,
were very poorly conditioned).
Wilkinson writes~\cite[p.~3]{Wilkinson:84}: ``Speaking for myself
I regard it as the most traumatic experience in my career as a
numerical analyst.''

Leaving aside the issue of why it was not expected that the
perfidious polynomial would be ill-conditioned, it
should suffice to us to observe
that the phenomenon of ill-conditioned characteristical polynomials
was so common that the practice of computing eigenvalues
via finding their zeros was completely abandoned on the grounds
of this ill-conditioning. This is explicit in the following passage
by Trefethen and Bau~\cite[p.~190]{TrefethenBau}:

\begin{quote}
{\small
Perhaps the first method one might think of would be to compute
the coefficients of the characteristic polynomial and use a rootfinder
to extract its roots. Unfortunately [\dots] this strategy is a bad one,
because polynomial rootfinding is an ill-conditioned problem in
general, even when the underlying eigenvalue problem is
well-conditioned.
}
\end{quote}
While this passage is right to point at numerical stability as
the stumbling block for the computation of eigenvalues
via root-finding for characteristic polynomials, its choice of
words is somehow hapless as the statement ``polynomial 
rootfinding is an ill-conditioned problem in general'' may lead 
to the impression that there is a natural choice for all the 
involved ingredients, notably, for the basis we take for 
the space of degree $n$ polynomials, for the probability 
distribution we draw the resulting coefficients from, and  
for the measure of condition used.

Maybe the least controversial of these ingredients is 
the probability distribution, as it is common practice to use 
the standard Gaussian (this goes back to the origins
of modern numerical analysis~\cite{vNGo51} and is 
more recently found in the work by Borgwart~\cite{Borgwardt4},
Demmel~\cite{Demmel88}, and Smale~\cite{Smale97}
among others). 
The most common choice for condition measures 
is the {\em normwise relative} version of the condition number. 
Coupled with some choice of bases (notably, the standard monomial 
basis) this condition number is known 
to be small in general (see, e.g.,~\cite[\S17.8]{Condition}). The 
fact that the last extends to multivariate polynomial systems 
is precisely what has allowed the recent
advances in zero-finding for these 
systems~\cite{BePa:11,BuCu11}. But there may as well be
bases with respect to which this is no longer the case. 
These facts are for the 
normwise relative measure of condition. It is a must to observe, 
at this point, that for another common measure of condition, the 
so called {\em componentwise}, they don't need to be true. We will 
return to this issue on \S\ref{sec:sim}. 

Returning to the eigenvalue problem, 
it is known that Gaussian matrices are, on the average,
well-conditioned with respect to the eigenpair
problem~\cite[Thm.~2.14]{ABBCS:15}.
And, as we just mentioned, typical polynomials are 
well-conditioned for the computation of their zeros 
w.r.t.\ the monomial basis and normwise relative condition 
number. 

Our goal in this paper is to show that in this setting {\em characteristic 
polynomials of typical matrices are ill-conditioned for the computation of 
their zeros}.

In order to precisely state the result, we recall that
the condition of a complex polynomial $f$
at its zero $\z\in\C$ is defined as
the worst possible magnification of the error
in the returned zero $\tilde\z$ with respect to the
size of a perturbation $\tilde{f}$ of $f$. In the  
relative normwise setting these errors are measured
normwise for the polynomial $f$ and they are relative
for both $f$ and $\z$ (see e.g.,~\cite[\S3]{Demmel:87b}). 
That is,
\begin{equation}\label{eq:def-condition}
   \cond(f,\z):=\lim_{\delta\to0} \sup_{\|\tilde{f}-f\|\leq \delta}
   \frac{|\tilde\z-\z|}{\|\tilde{f}-f\|}\,\frac{\|f\|}{|\z|}.
\end{equation}

In what follows we assume that the norm above is 
with respect to the standard monomial basis.  
We also denote by $\chi_A$ the characteristic
polynomial of a complex matrix $A\in\Cnn$, and
we say that $A$ is {\em standard Gaussian}
when the real and imaginary parts of its entries are independent
standard Gaussian random variables. A result by
Kostlan (Theorem~\ref{thm:kostlan} below) shows that for
such a matrix we can individualize its $n$ different eigenvalues
by the distribution of their moduli.

The main result of this paper is the following.

\begin{theorem}\label{th:mainR}
Suppose that $A\in\Cnn$ is standard Gaussian and
$\lambda_1,\ldots,\lambda_n$ are the eigenvalues of $A$.
For $1\le i\le n$ we have
$$
 \E \ln \cond(\chi_A, \lambda_i) \ \geq\
  \frac12 (n-1) \ln i - 0.79\, n - 0.5\, i .
$$
Moreover, there exists $K>0$ such that for all~$n$ we have
$$
 \min_{1\le i \le n} \E \ln \cond(\chi_A, \lambda_i) \ \geq\ 0.05\, n  -  K .
$$
For the average logarithm of the condition number we obtain
$$
 \E\Big(\frac{1}{n} \sum_{i=1}^n \ln\cond(\chi_A, \lambda_i) \Big)
 \ \geq\  \frac12 (n-1)\ln n - 1.54\,n .
$$
\end{theorem}

The third bound in Theorem~\ref{th:mainR} can be interpreted
in terms of the so called {\em standard distribution} on the solution
manifold. The latter is the set
$$
   V:=\{(A,\lambda)\in\Cnn\times \C\mid \det(A-\lambda\Id)=0\}
$$
and the standard distribution on $V$ amounts to drawing
$A$ from the standard Gaussian and then drawing one of its
(almost surely) $n$ different eigenvalues from the uniform
distribution. We denote this standard distribution on $V$
by $\rho_{\std}$. The third bound can then be written as
$$
  \E_{(A,\lambda)\sim\rho{\std}}\ln\cond(A,\lambda)
  \geq \frac12 (n-1)\ln n - 1.54\,n.
$$
Note that this implies, via Jensen's inequality, that
$$
  \E_{(A,\lambda)\sim\rho{\std}}\cond^2(A,\lambda)
  \geq n^{n-1} e^{-3.08\,n}.
$$
The following result is a small improvement on this lower
bound.

\begin{theorem}\label{thm:average-standard}
$$
  \E_{(A,\lambda)\sim \rho_{\std}}
  (\cond^2(\chi_A, \lambda)) \geq (n-1)!\,2^n.
$$
\end{theorem}

After laying down some preliminaries, we
prove these results in Section~\ref{se:results}.
Then, in Section~\ref{sec:simul}, we discuss the robustness
of Theorem~\ref{th:mainR} with respect to (some) changes
in the way errors are measured for conditioning, and
show some computer simulations. The latter are consistent with
our findings.
\medskip

\noindent
{\bf Acknowledgment.\quad}
We are grateful to Dennis Amelunxen for performing the
the computer simulations.

\section{Preliminaries}
\label{se:prelim}

\subsection{Condition of univariate polynomials}
\label{se:cond-univar}

The Euclidean norm of a complex polynomial $f\in\C[X]$,
\begin{equation}\label{eq:f}
      f(X)= a_nX^n+a_{n-1}X^{n-1}+\cdots+a_1X+a_0,
\end{equation}
is given by
\begin{equation*}\label{eq:Enorm}
  \|f\|^2:= \sum_{k=0}^n |a_k|^2.
\end{equation*}
For a simple  zero $\z\in\C$ of $f$ the condition number
$\cond(f,\z)$ in~\eqref{eq:def-condition} then takes
the following form (see, e.g., in~\cite[\S14.1.1]{Condition})
\begin{equation}\label{eq:def-cond}
   \cond(f,\z)=\frac{\|f\|}{|\z|}\,\frac{1}{|f'(\z)|}\,
    \big\|(1,|\z|,\ldots,|\z|^n)\big\|.
\end{equation}

We won't use the following result in the proof of
Theorems~\ref{th:mainR} and~\ref{thm:average-standard}.
But as we have not seen it
anywhere and it has an interest per se, we offer here a proof.

\begin{proposition}\label{prop:lb-cond}
If $f$ is a complex polynomial of degree $n$ and $\z\in\C$ a zero
of $f$, then $\cond(f,\z) \ge \frac{1}{n}$.
\end{proposition}

\proof
Write $r=|\z|$ so that we have
$$
   \cond(f,\z)=\frac{\|f\|}{|f'(\z)|}\cdot
   \frac{\|(1,r,r^2,\ldots,r^n)\|}{r}.
$$
Assume $f$ is as in~\eqref{eq:f}. Then
\begin{eqnarray*}
  |f'(\z)| &=& \big|na_n\z^{n-1} + \ldots+ 2a_2\z+a_1\big| \\
 &\leq& \big\|(a_n,a_{n-1},\ldots,a_1)\big\|\cdot
  \big\|(n\z^{n-1},(n-1)\z^{n-2},\ldots,2\z,1)\big\|\\
&\leq & \|f\| \big\|(nr^{n-1},(n-1)r^{n-2},\ldots,2r,1)\big\|
 \,\leq\, \|f\|\, n\, \big\|(r^{n-1},r^{n-2},\ldots,r,1)\big\|,
\end{eqnarray*}
the first line by Cauchy-Schwartz. It follows that
\begin{equation}\tag*{\qed}
\cond(f,\z)\;\geq\; \frac{\|(1,r,r^2,\ldots,r^n)\|}
{nr\,\|(r^{n-1},r^{n-2},\ldots,r,1)\|}
\;\geq\; \frac{\|(1,r,r^2,\ldots,r^n)\|}
{n\,\|(r,r^2,\ldots,r^n)\|} \;\geq\; \frac1n.
\end{equation}

\subsection{Distribution of eigenvalues of Gaussian matrices}

Let $A=(z_{ij})\in\Cnn$ be a random matrix such that for all $i,j$,
the real part $\Re(z_{ij})$ and the imaginary  part $\Im(z_{ij})$
of $z_{ij}$ are independent standard Gaussian random variables.
Note that $\E |z_{ij}|^2 = 2$.
The resulting distribution of matrices is sometimes called the
{\em complex Ginibre ensemble}.
We will also say that $A$ is {\em standard Gaussian}
and write $A\sim N(0,\Id)$ for this.
Ginibre~\cite{Ginibre} showed that the density of the
joint probability distribution of the eigenvalues
$\lambda_1,\ldots,\lambda_n$ of $A$ is given by
\begin{equation}\label{eq:ginibre}
  \rho(\lambda_1,\ldots,\lambda_n)
  = C_n\, e^{-\frac12\sum_{i=1}^n |\lambda_i|^2}
  \prod_{i<j} |\lambda_i - \lambda_j|^2 ,
\end{equation}
where $C_n^{-1} = 2^{\frac{n(n+1)}{2}} \pi^n \prod_{k=1}^{n} k! $.

Based on Ginibre's formula, Eric Kostlan~\cite{Kostlan92} observed
that, surprisingly, the squared absolute values~$|\lambda_i|^2$ of the
eigenvalues~$\lambda_i$ of a standard Gaussian matrix $A$ are
distributed like independent $\chi^2$ random variables.  This
insight will be crucial for our analysis.

\begin{theorem}\label{thm:kostlan}
For a standard Gaussian $A\in\C^{n\times n}$ with eigenvalues $\lambda_1,\ldots,\lambda_n$,
the set $\{|\lambda_1|^2,\dots,|\lambda_n|^2\}$
is distributed like the set $\{\chi_{2}^2, \dots, \chi_{2n}^2\}$,
where $\chi^2_{2},\ldots,\chi^2_{2n} $ are
independent $\chi^2$ random variables with 
$2,\ldots,2n$ degrees of freedom,
respectively.
\eproof
\end{theorem}

\subsection{Some useful bounds}

We collect here some known facts related to $\chi^2$-distributions
needed for the proof of the main result.

The {\em psi function}, also called the {\em logarithmic digamma function}, is defined as
the logarithmic derivative of the gamma function:
$$
 \psi(x) := \frac{d}{dx} \ln\Gamma(x) = \frac{\Gamma'(x)}{\Gamma(x)} .
$$
It satisfies the recursion
$\psi(x+1) = \psi(x) + \frac{1}{x}$ for $x>0$, and
$\psi(1) = -\gamma \approx -0.577$, where $\gamma$ denotes
the Euler-Mascheroni constant.
Therefore, for positive $m,n\in\N$,
\begin{equation}\label{eq:psi-values}
\psi(n) =  - \gamma + \sum_{k=1}^{n-1} \frac{1}{k} \ \ge\ \ln n - \gamma
\qquad\mbox{and}\qquad
\psi(m+n) - \psi(m) = \sum_{k=m}^{m+n-1} \frac{1}{k} .
\end{equation}

We say that  a nonnegative random variable $X$ is $\chi$-distributed with $n$~degrees of freedom,
written $X\sim \chi_n$, if $X^2\sim \chi^2_n$.

\begin{lemma}\label{le:error2zero}
Suppose that $r_1,r_2,\ldots$ is a sequence of independent $\chi$-distributed random variables,
where $r_i\sim\chi_{2i}$. Then we have for $i,j\ge 1$
$$
 \E\ln r_i^2 \ =\ \psi(i) + \ln 2 \ \ge\
   \E\ln r_1^2 = \ln 2 - \gamma >  0.1159, 
$$
and hence
$$
 \E\ln \frac{r_j^2}{r_i^2 + r_j^2} = \psi(j) - \psi(i+j).
$$
Moreover, if $j\ge 2$, we have
$$
   \E\ln\frac{r_i + r_j}{r_j} \ \le\ \sqrt{\frac{i}{j-1}} .
$$
\end{lemma}

\proof
The density of $\chi^2_n$ is given by
$\rho(q) = 2^{-\frac{n}{2}}\Gamma(\frac{n}{2})^{-1} q^{\frac{n}{2}-1}
e^{-\frac{q}{2}}$.
Therefore, substituting $v=q/2$,
\begin{equation*}
\begin{split}
 \E \ln\chi^2_n &= \frac{1}{2^{\frac{n}{2}}\Gamma(\frac{n}{2})}
   \int_0^\infty q^{\frac{n}{2}-1} e^{-\frac{q}{2}} \ln(q)\, dq \ = \
  \frac{1}{\Gamma(\frac{n}{2})}
   \int_0^\infty v^{\frac{n}{2}-1} e^{-v} (\ln(v) + \ln(2)) \, dq \\
  & = \psi(n/2) + \ln(2) ,
\end{split}
\end{equation*}
where we used~\cite[4.352-1]{GradRy} for the last equality.
This shows the first assertion.

For the second assertion we use
$\ln \frac{r_i + r_j}{r_j} = \ln\big(1+\frac{r_i}{r_j}\big) \le \frac{r_i}{r_j}$.
Note that $\frac{r^2_i/i}{r^2_j/j}$ has the law of an
F-distribution with parameters $2i$ and $2j$,
whose expectation is
known~\cite[p.~326]{JKB2:95}
to be equal to $2j/(2j-2)$. Therefore
\begin{equation}\tag*{\qed}
 \E \frac{r_i}{r_j} \ \le\ \sqrt{\E \frac{r^2_i}{r^2_j}} \ =\
  \sqrt{\frac{i}{j-1}} .
\end{equation}

\begin{lemma}\label{le:novo}
Suppose that $r_1,r_2,\ldots$ is a sequence of independent $\chi$-distributed random variables,
where $r_i\sim \chi_{2i}$. Then we have for $i\ge 2$ and $j\ge 1$
$$
  \E\ln \frac{r_i r_j}{r_i + r_j} \ \ge\ \E\ln \frac{r_2 r_j}{r_2 + r_j}  - \frac12\ln 2 .
$$
\end{lemma}

\proof
We first note that for $i\ge 2$ and $j\ge 1$:
\begin{equation}\label{eq:psi}
 \psi(i) +\psi(j) - \psi(i+j) \ge \psi(2) + \psi(j) - \psi(2+j) .
\end{equation}
Indeed, this means
$\psi(i+j)-\psi(i)\leq \psi(2+j)-\psi(2)$, which
by \eqref{eq:psi-values} is equivalent to
$$
  \sum_{k=i}^{i+j-1}\frac{1}{k} \ =\
  \sum_{\ell=2}^{j+1} \frac{1}{\ell+(i-2)}  \ \le\
  \sum_{\ell=2}^{j+1}\frac{1}{\ell},
$$
which is obviously true.

Using Lemma~\ref{le:error2zero} and inequality~\eqref{eq:psi}
we deduce that
\begin{equation}\label{eq:EF}
  \E\ln \frac{r_i^2r_j^2}{r_i^2 + r_j^2} \ \ge\
  \E\ln \frac{r_2^2r_j^2}{r_2^2 + r_j^2}.
\end{equation}
Now we use the fact
$r_i^2 +r_j^2 \le (r_i +r_j)^2 \le 2\, (r_i^2 +r_j^2)$
to obtain, for $i\ge 2$,
\begin{align}
 \E\ln\frac{r_ir_j}{r_i + r_j}  &=\;
 \frac12 \E\ln \frac{r_i^2r_j^2}{(r_i+r_j)^2}
 \;\ge\; \frac12 \E\ln \frac{r_i^2r_j^2}{r_i^2+r_j^2} -\frac12\ln 2
 \nonumber\\
 &\stackrel{\eqref{eq:EF}}{\ge}\;
 \frac12 \E\ln \frac{r_2^2r_j^2}{r_2^2+r_j^2} -\frac12\ln 2
 \;\ge\;
 \frac12 \E\ln \frac{r_2^2r_j^2}{(r_2+r_j)^2} -\frac12\ln 2
 \nonumber\\
 &=\; \E\ln \frac{r_2r_j}{r_2+r_j} -\frac12\ln 2.
     \tag*{\qed}
\end{align}

\section{Proofs of the Main Results}
\label{se:results}

\proofof{Theorem~\ref{th:mainR}}
We denote by $\chi_A$ the characteristic polynomial of a complex
matrix $A\in\Cnn$ and by $\{\lambda_1,\ldots,\lambda_n\}$ the
multiset
of its eigenvalues, so that $\chi_A = (X-\lambda_1)\cdots
(X-\lambda_n)$.  We rely on Theorem~\ref{thm:kostlan} to
associate the index $i$ to the eigenvalue
satisfying $|\lambda_i|^2 \sim \chi_{2i}^2$.

In a first step we provide a lower bound for the condition of the pair
$(\chi_A,\lambda_i)$ as characterized in \eqref{eq:def-cond}.  Since
$(-1)^n\det(A)$ is the constant term of $\chi_A$, we have
$$
   \|\chi_A\| \ \ge\ |\det A | = |\lambda_1|\cdots |\lambda_n|.
$$
We also use the facts that, for each $1\le i\le n$, we have
$\big\|(1,|\lambda_i|,\ldots,|\lambda_i|^n)\big\| \ge
|\lambda_i|^{n-1}$ and
$$
   |\chi_A'(\lambda_i)|=\prod_{j\neq i}|\lambda_j-\lambda_i| .
$$
Replacing the last three relations in~\eqref{eq:def-cond}
we obtain, for each $1\le i\le n$,
\begin{equation}\label{eq:condchiA}
\begin{split}
   \cond(\chi_A,\lambda_i) &\ =\ \frac{\|\chi_A\|}{|\lambda_i|}\,
    \frac{1}{|\chi'_A(\lambda_i)|}\,
    \big\|(1,|\lambda_i|,\ldots,|\lambda_i|^n)\big\| \\
  & \ \ge\ \prod_{j\ne i} \, \frac{|\lambda_i||\lambda_j|}
     {|\lambda_i-\lambda_j|}
   \ \ge\ \prod_{j\ne i} \, \frac{|\lambda_i||\lambda_j|}
    {|\lambda_i| + |\lambda_j|} ,
\end{split}
\end{equation}
the last by the triangle inequality.

In what follows we write $r_j:=|\lambda_j|$, for $j=1,\ldots,n$.
Recall, these are
independent random variables with $r_j\sim\chi_{2j}$.
From~\eqref{eq:condchiA} we get for fixed $i\ge 1$,
$$
   \ln \cond(\chi_A,\lambda_i) \ \ge\ \sum_{j\ne i}
    \ln\frac{r_ir_j}{r_i + r_j} ,
$$
and hence
\begin{equation}\label{eq:ELB}
  \E\ln \cond(\chi_A,\lambda_i) \ \ge\ \sum_{j\ne i}
 \Big(\E\ln r_i + \E\ln \frac{r_j}{r_i + r_j} \Big) .
\end{equation}
To bound the first term in the right-hand side,
we combine Lemma~\ref{le:error2zero}
with~\eqref{eq:psi-values} to obtain
\begin{equation}\label{eq:suppo}
 \E\ln r_i = \frac12 \psi(i) + \frac12 \ln 2 \ \ge\
  \frac12 \ln i -\frac12\gamma + \frac12 \ln 2.
\end{equation}
We bound the second term in the right-hand side
of~\eqref{eq:ELB} using that
$(r_i +r_j)^2 \le 2(r_i^2 +r_j^2)$,
$$
 \ln \frac{r_j}{r_i + r_j} \ =\
 \frac12\ln \frac{r_j^2}{(r_i + r_j)^2} \ \ge\
 \frac12 \ln\frac{r_j^2}{2(r_i^2 + r_j^2)}
 = \frac12 \ln\frac{r_j^2}{r_i^2 + r_j^2} - \frac12 \ln 2 .
$$
Using Lemma~\ref{le:error2zero} and~\eqref{eq:psi-values}
one more time, we deduce
$$
 \E \ln \frac{r_j}{r_i + r_j} \ \ge\
   \frac12  \big(\psi(j) - \psi(i+j) \big)  - \frac12 \ln 2 \ =\
  - \frac12 \sum_{k=j}^{i+j-1}\frac{1}{k} - \frac12 \ln 2 .
$$
Replacing this bound and~\eqref{eq:suppo}
in~\eqref{eq:ELB}, we obtain
\begin{equation}\label{eq:Star}
 \E\ln \cond(\chi_A,\lambda_i) \ \ge\
 \frac12 (n-1) \ln i - \frac{\gamma}{2} (n-1)
  - \frac12 \sum_{j=1}^n \sum_{k=j}^{i+j-1} \frac{1}{k} .
\end{equation}
We can further bound the last term in the right-hand term
by exchanging the order of summation.
For fixed $1\le k \le i + n- 1$ we have at most $k$ choices of $j$,
since $1\le j \le k$. Therefore,
$$
   \frac12 \sum_{j=1}^n \sum_{k=j}^{i+j-1} \frac{1}{k}  \ \le \
   \frac12 \sum_{k=1}^{n+i-1} \frac{k}{k}  = \frac12 (n+i-1).
$$
We thus obtain from~\eqref{eq:Star}
\begin{equation}\label{eq:NU}
\begin{split}
   \E\ln \cond(\chi_A,\lambda_i) &\ge\ \frac12 (n-1) \ln i
   -\frac{\gamma}{2} n + \frac{\gamma}{2} - \frac{n+i-1}{2} \\
  & =\ \frac12 (n-1) \ln i - \frac{1+\gamma}{2} n
     - \frac{i}{2} + \frac{1+\gamma}{2} \\
  & \ge\ \frac12 (n-1) \ln i - 0.79\, n - 0.5\, i ,
\end{split}
\end{equation}
which proves the first assertion.
Averaging this over $i=1,\ldots,n$, we obtain
\begin{equation*}
\begin{split}
 \frac{1}{n} \sum_{i=1}^n \E\ln \cond(\chi_A,\lambda_i) & \ge\
 \frac12 (n-1) \frac{1}{n } \sum_{i=1}^n \ln i - 0.79\, n - 0.5\,
    \frac{1}{n}\sum_{i=1}^n i \\
  & \ge\ \frac12 (n-1) (\ln n -1)  - 0.79\, n - \frac{n+1}{4},
\end{split}
\end{equation*}
where we have used that
$$
 \sum_{i=1}^n \ln i \ \ge\ \int_1^n \ln x\, dx \ \ge\  n(\ln n -1).
$$
Therefore,
$$
  \frac{1}{n} \sum_{i=1}^n \E\ln \cond(\chi_A,\lambda_i)  \ \ge\
  \frac12 (n-1)\ln n  - 1.54\, n,
$$
which proves the third assertion.

Note that the general lower bound \eqref{eq:NU} is not useful
for small values of~$i$. This is why the proof of
the second assertion on the minimum of the
expectations over $1\le i \le n$ needs an extra argument.

We first consider the case where $i=1$.
Lemma~\ref{le:error2zero} tells us that $\E \ln r_1  > 0.057$ and
$$
0 \ \le\ \E\ln \frac{r_1 + r_j}{r_j} \ \le\ \sqrt{\frac{1}{j-1}} \ <\ 0.007 ,
$$
provided $j\ge j_1:= 20410$, which implies in that case
$$
 \E\ln \frac{r_1r_j}{r_1+ r_j} = \E \ln r_1 + \E\ln\frac{r_j} {r_1 + r_j}
  \ \ge\ 0.057 - 0.007 = 0.05 .
$$
If we denote by $K_1$ the sum of $\E\ln \frac{r_1 + r_j}{r_1r_j}$
over $1\le j \le j_1-1$, then we obtain with~\eqref{eq:ELB} that,
for all $n\ge j_1$,
\begin{equation}\label{eq:uno}
  \E \ln \cond(\chi_A,\lambda_1) \ \ge\  -K_1 +
  \sum_{j=j_1}^n 0.05 \ = 0.05 \cdot n - (K_1 + 0.05\, j_1) .
\end{equation}
We study now the case where $2\le i\le n$.
Combining~\eqref{eq:ELB} with Lemma~\ref {le:novo} we get
$$
  \E\ln \cond(\chi_A,\lambda_i) \ \ge\ \sum_{j\ne i}
  \E\ln \frac{r_ir_j}{r_i + r_j} \ \ge\
  \sum_{j\ne i} \Big(\E\ln \frac{r_2r_j}{r_2 + r_j} -\frac12 \ln 2 \Big) .
$$
We have
$\E\ln r_2 -\frac12\ln 2 = \frac{1-\gamma}{2} >0.2$
by Lemma~\ref{le:error2zero}.
Choose $j_2 = 180$ such that
$\sqrt{\frac{2}{j_2-1}} \le \frac{1-\gamma}{4}$.
Using again Lemma~\ref{le:error2zero} we obtain,
for $j\ge j_2$,
$$
 \E\ln \frac{r_2r_j}{r_2 + r_j} -\frac12\ln 2
  = \E\ln r_2 -\frac12\ln 2 + \E\ln \frac{r_j}{r_2 + r_j}
\ge \frac{1-\gamma}{2} - \sqrt{\frac{2}{j-1}} \ge \frac{1-\gamma}{4} .
$$
Writing
$$
 K_2 := - \min_{2\le i < j_2}
  \sum_{1\le j< j_2 \atop j\ne i} \Big(\E\ln\frac{r_2r_j}{r_2+r_j}
  - \frac12\ln 2 \Big)
$$
we conclude that for all $2\le i \le n$ and all $n\ge j_2$,
\begin{equation*}\label{eq:due}
 \E \ln \cond(\chi_A,\lambda_i)\ge
-K_2 + \sum_{j_2\le j\le n \atop j\ne i} \frac{1-\gamma}{4}
 \ge \frac{1-\gamma}{4} (n -j_2) - K_2
 \geq 0.1\, n-(0.11\, j_2+K_2).
\end{equation*}
The second assertion of the theorem follows from this bound
and~\eqref{eq:uno} taking
$K:=\max\{0.05\,j_1+K_1,0.11\, j_2+K_2\}$.
\eproof

\proofof{Theorem~\ref{thm:average-standard}}
The joint density of the eigenvalues of a Gaussian matrix
$A$ given in~\eqref{eq:ginibre} is invariant under permutations
of the $\lambda_i$s. Hence, for all $i\leq n$,
$$
  \E(\cond^2(\chi_A, \lambda_i))=\E(\cond^2(\chi_A, \lambda_n)).
$$
We will therefore compute the expectation for $i=n$ (note that
this is just for notational convenience: $\lambda_n$ now is not
the eigenvalue whose modulus squared is
$\chi^2_{2n}$-distributed).

We know that
$$
  \cond^2(\chi_A, \lambda_n)\geq
 \frac{\prod_{j=1}^n |\lambda_j|^2}{|\lambda_n|^2}\;
 \frac{|\lambda_n|^{2n}}{\prod_{j<n} |\lambda_n-\lambda_j|^2}
 =
 \prod_{j<n} |\lambda_j|^2\; \frac{|\lambda_n|^{2n}}{\prod_{j<n}
 |\lambda_n-\lambda_j|^2}.
$$
Then, because the density of $(\lambda_1, \dots, \lambda_n)$ is (\ref{eq:ginibre}),
\begin{eqnarray*}
& &\E(\cond^2(\chi_A, \lambda_n))\\
&\geq&
 \int \prod_{j<n} |\lambda_j|^2 \frac{|\lambda_n|^{2n}}{\prod_{j<n}
 |\lambda_n-\lambda_j|^2}\, C_n\,
e^{-\frac12\sum_{i=1}^n |\lambda_i|^2}\prod_{j<i}
|\lambda_i - \lambda_j|^2 \,d\lambda_1\dots d\lambda_n\\
&=&
  \int C_n e^{-\frac12\sum_{i=1}^{n-1} |\lambda_i|^2}
 \prod_{j<n} |\lambda_j|^2  \prod_{j<i< n}
 |\lambda_i - \lambda_j|^2 \,d\lambda_1\dots d\lambda_{n-1}
 \cdot \int  |\lambda_n|^{2n} e^{-\frac12|\lambda_n|^2}  \, d\lambda_n.
\end{eqnarray*}
The first integral yields
$$
  \frac{C_n}{C_{n-1}} \E_{A\sim N(0, \Id_{n-1})}
  (|\lambda_1|^2\cdots|\lambda_{n-1}|^2) =
  \frac{C_n}{C_{n-1}} 2^{n-1}(n-1)! =
  \frac{1}{2\pi n}
$$
because of Theorem~\ref{thm:kostlan} and the fact
that $\frac{C_n}{C_{n-1}}=\frac{1}{\pi n!2^n}$. 

The second integral yields
\begin{eqnarray*}
  \int  |\lambda|^{2n} e^{-\frac12|\lambda|^2}  \, d\lambda 
   &=& 2\pi \int_0^\infty r^{2n+1}e^{-\frac12 r^2}dr \\
   &=& 2\pi \int_0^\infty (2u)^{n}e^{-u}du = \pi 2^{n+1} 
   \Gamma(n+1)= \pi 2^{n+1}n!, 
\end{eqnarray*}
where we changed variables $r=|\lambda|$ and $u=r^2/2$ 
in the first and second equalities, respectively.
The result follows.
\eproof

\section{Numerical simulations and additional remarks}
\label{sec:simul}

\subsection{Numerical simulations}\label{sec:sim}

We have performed some
computer experiments to gauge the actual behavior of conditioning
for the characteristic polynomials of typical matrices.
Specifically, for each $n=2,\ldots,100$, we have generated
10,000 Gaussian matrices in $\Cnn$ and computed, in each case,
the average of the logarithms of the following two quantities:
$$
  \condmin(\chi_A):=\min_{i\leq n} \cond(\chi_A,\lambda_i)
$$
and
$$
  \condmax(\chi_A):=\max_{i\leq n} \cond(\chi_A,\lambda_i).
$$
Then, to visualize the growth of these averages with $n$
we have plotted (smoothed curves corresponding to) the
graphs of the following functions of $n$,
$$
  \frac{\Avg (\ln\condmin(\chi_A))}{n}
  \qquad\hbox{and}\qquad
  \frac{\Avg (\ln\condmax(\chi_A))}{n}.
$$
The resulting figure looks as follows.
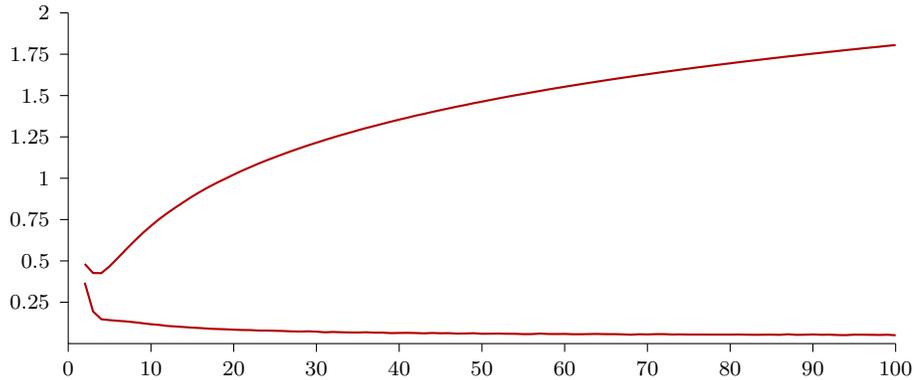
\begin{center}
\def\myXsc{0.05}
\begin{tikzpicture}[scale=2.2,xscale=\myXsc]
  \def\myE{0.05}
  \draw (0,0) -- (100,0) (0,0) -- (0,2);
  \foreach \n in {0,10,...,100}
    \draw (\n,0) -- (\n,-\myE) node[below]{$\scriptstyle\n$};
  \foreach \y in {0.25,0.5,...,2}
    \draw (0,\y) -- (-\myE/\myXsc,\y) node[left]{$\scriptstyle\y$};
  \draw[thick,red] plot file{fmax.table};
  \draw[thick,red] plot file{fmin.table};
\end{tikzpicture}
\end{center}
We observe that, on the average, $\frac{\ln\condmin(\chi_A)}{n}$
stabilizes on a value around $0.05$ whereas
$\frac{\ln\condmax(\chi_A)}{n}$ grows in what appears to be a
logarithmic manner. We therefore plotted the curve for
$\frac{\Avg (\ln\condmax(\chi_A))}{n\ln(n)}$ and still observed a
very gentle growth. So we finally did so for
$\frac{\Avg (\ln\condmax(\chi_A))}{n\ln(n)\ln(\ln n)}$
(and $n\geq 4$) and obtained the following figure.
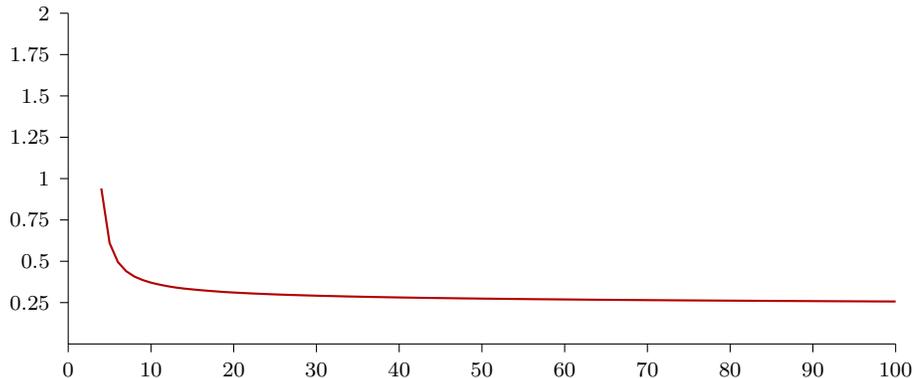
\begin{center}
\def\myXsc{0.05}
\begin{tikzpicture}[scale=2.2,xscale=\myXsc]
  \def\myE{0.05}
  \draw (0,0) -- (100,0) (0,0) -- (0,2);
  \foreach \n in {0,10,...,100}
    \draw (\n,0) -- (\n,-\myE) node[below]{$\scriptstyle\n$};
  \foreach \y in {0.25,0.5,...,2}
    \draw (0,\y) -- (-\myE/\myXsc,\y) node[left]{$\scriptstyle\y$};
  \draw[thick,red] plot file{fmaxloglog.table};
\end{tikzpicture}
\end{center}
The value of $\frac{\Avg (\ln\condmax(A))}{n\ln(n)\ln(\ln n)}$
appears to stabilize at around $0.25$ but we note that the
numerical results here are not enough to determine whether
$\E\log\condmax(\chi_A)$ grows as
$n\ln(n)$ or as $n\ln(n) \ln(\ln n)$.
\medskip

The condition number we have considered in all the previous
development is defined in a {\em relative normwise} manner.
It measures errors in the approximation $\tilde{f}$ of a
polynomial $f$ by the  quotient $\frac{\|\tilde{f}-f\|}{\|f\|}$.
It should come as no surprise that in the case of $\chi_A$
with $A$ Gaussian, the condition number $\cond(\chi_A,\lambda)$
will be large for each~$\lambda$. After all, we are allowing
errors proportional to $|\det(A)|$ in coefficients which we expect
to be much smaller than this determinant. A different way to
measure the error in $f$ is the {\em componentwise}. If
$f=a_n X^n+\cdots+a_1X+a_0$ then we measure the error
of an approximation $\tilde{f}$ by
$\max_{0\leq i\leq n} \frac{|\tilde{a}_i-a_i|}{|a_i|}$. This
leads to the componentwise condition number
$$
   \Cw(f,\zeta):=\lim_{\delta\to0} \sup_{|\tilde{a_i}-a_i|\leq \delta\,|a_i|}
   \max_{0\leq i\leq n}\frac{|\tilde\z-\z|}{|\tilde{a_i}-a_i|}\,\frac{|a_i|}{|\z|}.
$$
It turns out that this condition number has a simple
characterization (take $m=1$
in~\cite[Thm.~2.1]{Winkler:06})
$$
  \Cw(f,\z)=\frac{1}{|\z|}\,\frac{1}{|f'(\z)|}\,
     \sum_{i=0}^n |a_i|\,|\z|^i.
$$
We have used this expression to perform some computations,
similar to the preceding ones (10,000 Gaussian matrices
in $\Cnn$ for each of $n=2,\ldots,100$), to gauge the average
behavior of
$$
  \Cw_{\max}(\chi_A):=\max_{i\leq n}\Cw(\chi_A,\lambda_i).
$$
The following picture, plotting
$\frac{\Avg(\ln\Cw_{\max}(\chi_A))}{\ln n}$,
suggests that $\Cw_{\max}(\chi_A)$ has a sublinear
(or maybe linear) growth.
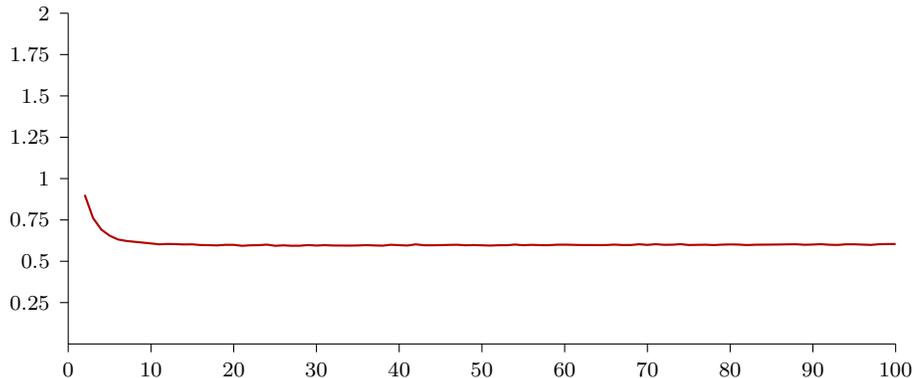
\begin{center}
\def\myXsc{0.05}
\begin{tikzpicture}[scale=2.2,xscale=\myXsc]
  \def\myE{0.05}
  \draw (0,0) -- (100,0) (0,0) -- (0,2);
  \foreach \n in {0,10,...,100}
    \draw (\n,0) -- (\n,-\myE) node[below]{$\scriptstyle\n$};
  \foreach \y in {0.25,0.5,...,2}
    \draw (0,\y) -- (-\myE/\myXsc,\y) node[left]{$\scriptstyle\y$};
  \draw[thick,red] plot file{fmaxcomploghalf.table};
\end{tikzpicture}
\end{center}
The good behavior of $\Cw_{\max}(\chi_A)$ stands in contrast
with the numerical instability observed when computing the zeros
of characteristic polynomials. The most likely explanation is that
the algorithms that compute $\chi_A$ from $A$ produce
forward-errors that are {\em not} componentwise small.
We have browsed the literature in search of some
bound for this forward-error (for some algorithm) and have
found none. The closest result we found is in a paper by
Ipsen and Rehman that gives upper
bounds (Theorem~3.3 and Remark~3.4 in~\cite{IpRe:08})
for componentwise errors in $\chi_A$
due to normwise measured perturbations on $A$.
These upper bounds are absolute
(i.e., not relative to the moduli of the coefficients
of $\chi_A$). The corresponding relative bounds
are small for the ``extreme'' coefficients (corresponding to
large and small degree) but may be large for the ``middle''
coefficients. But it is unclear whether this is so because the
relative condition number for these coefficients
are large or because the upper
bounds in~\cite{IpRe:08} are not sharp.

\subsection{Additional remarks}\label{sec:additional}

Condition numbers depend on the way errors (for both input data
and output) are measured. When norms are used to do so, the
choice of a particular norm plays a role as well. In all our previous
development we have used the Euclidean norm in the space
of polynomials. For reasons to be explained soon, another
common choice in this space is the Weyl norm, which, for
$f$ as in~\eqref{eq:f}, is given by
\begin{equation}\label{eq:Wnorm}
  \|f\|_W^2:= \sum_{k=0}^n {{n}\choose{k}}^{-1}|a_k|^2.
\end{equation}
Again, for a simple  zero $\z\in\C$ of $f$, the condition number
$\cond_W(f,\z)$ induced by the Weyl norm
is shown (see, e.g.,~\cite[\S14.1.1]{Condition}) to have the form
\begin{equation}\label{eq:def-condW}
   \cond_W(f,\z)=\frac{\|f\|_W}{|\z|}\,\frac{1}{|f'(\z)|}\,
    \big(1+|\z|^2\big)^{\frac{n}2}.
\end{equation}
Since $\|\chi_A\|_W\geq |\det(A)|$ and
$\big(1+|\lambda_i|^2\big)^{\frac{n}2}\geq |\lambda_i|^{n-1}$
we have
$$
  \cond_W(\chi_A,\lambda_i)\geq \prod_{j\ne i} \,
  \frac{|\lambda_i||\lambda_j|}{|\lambda_i| + |\lambda_j|}.
$$
That is,~\eqref{eq:condchiA}  holds for
$\cond_W(\chi_A,\lambda_i)$ and, hence, Theorem~\ref{th:mainR}
holds as well (as the proof of this theorem is just a bound for
the quantity on the right-hand side above).

Yet another change of setting consists of homogenizing the
polynomial $f$ so that its zeros $\xi$ now live in $\P^1(\C)$. For
the homogenization
$$
  f^{\mathsf{h}}(X,Z)= a_nX^n+ a_{n-1}X^{n-1}Z +\cdots+
    a_1XZ^{n-1}+a_0Z^n
$$
of $f$ in~\eqref{eq:f} it is common to consider the Weyl norm
above because it is invariant under the action of the unitary
group. That is, for all unitary matrices $U\in\C^{2\times 2}$,
$\|f^{\mathsf{h}}(X,Z)\|=\|f^{\mathsf{h}}(U(X,Z))\|_W$. The major
difference in the notion of condition, however, now comes from the
fact that errors in the zero $\xi$ are measured with the Riemannian
distance in $\P^1(\C)$, i.e., with angles between lines in $\C^2$.
The condition thus obtained is therefore not affected by the
distortions produced by the zeros of~$f$ in $\C$ becoming large
(in modulus). This condition number was introduced by
Mike Shub and Steve Smale (they named it
$\mu(f^{\mathsf{h}},\xi)$) who,
in addition, showed the following
characterization~\cite{Bez3}\footnote{The writing in terms of
$\cond_W(f,\xi)$ is ours. Furthermore, we note that the formula
for $\mu$ in~\cite[p.~6]{Bez3} has a typo ---$\|u\|^{d-1}$ should be
$\|u\|^d$--- and that we have disregarded the ``normalization
factor'' $d^{1/2}$.},
$$
  \mu(f^{\mathsf{h}},[\z\colon\!1]) =
  \cond_W(f,\z)\frac{|\z|}{\sqrt{1+|\z|^2}}.
$$
Here $\z\in\C$ and $[\z\colon\!1]$ is its image under the
standard inclusion $\C\hookrightarrow\P^1(\C)$. Then,
\begin{eqnarray*}
\mu(\chi_A^{\mathsf{h}},[\lambda_i\colon\!1]) &=&
   \frac{\|\chi^{\mathsf{h}}_A\|_W}{|\lambda_i|}\,
   \frac{1}{|\chi_A'(\lambda_i)|}\,
    \big(1+|\lambda_i|^2\big)^{\frac{n}2}
    \frac{|\lambda_i|}{\sqrt{1+|\lambda_i|^2}} \\
&\geq &
    |\det(A)|\,\frac{1}{|\chi_A'(\lambda_i)|}\,
    \big(1+|\lambda_i|^2\big)^{\frac{n-1}2}
   \geq \prod_{j\ne i} \,
  \frac{|\lambda_i||\lambda_j|}{|\lambda_i| + |\lambda_j|}.
\end{eqnarray*}
That is,~\eqref{eq:condchiA}  also holds for
$\mu(\chi^{\mathsf{h}}_A,[\lambda_i\colon\!1])$ and with it,
Theorem~\ref{th:mainR}.

{\small
  \bibliography{../../../book/book}
}

\end{document}